\documentclass[eqnumis, pdf]{hjms}
\usepackage{amsfonts}
\usepackage{amsmath}
\usepackage{footnote}
\usepackage{multicol}
\usepackage{booktabs}
\usepackage[flushleft]{threeparttable}
\usepackage[font=small,labelfont=bf]{caption}

\begin{document}

\markboth{T. Ashale, G.Abebe, K. Venkateswarlu}{Geometry of AL-monoids}

\title{Geometry of Autometrized Lattice Ordered Monoids}

\author{Tekalign Regasa$^{1}$, Girum Aklilu\coraut $^2$, Kolluru Venkateswarlu$^{3}$}

\address{$^1$Department of Mathematics, Addis Ababa University, Ethiopia, Email:{tekalign.regasa@aau.edu.et\\
	$^2$Department of Mathematics, Addis Ababa University, Ethiopia, Email:{girum.aklilu@aau.edu.et}\\
	$^3$ Department of Computer Science and System Engineering,Visakhapatnam, Andhra University, India, Email: {drkvenkateswarlu@gmail.com}}}
\maketitle
\begin{abstract}
In this paper, we study the geometry of Autometrized lattice ordered monoid(short as AL-monoids) by introducing the concept of metric betweeness and its properties $t_{1},t_{2},$ B-linearity, D-linearity, lattice betweeness, B-linearity, and D-linearity, segments and equilateral triangles. It is proved that there do not exist equilateral triangles in AL-monoids. It is also proved that any AL-monoid is ptolemaic. This result subsumes all the geometry of commutative DRl-semigroups.
\end{abstract}
\subjclass:{47H09, 47H10}  


\setlength\columnsep{15cm} 
\section{Introduction} 
The study of spaces where distances are chosen from algebraic structures other than real or complex numbers has received attention on occasion.
Since the metric operation $'$*$'$ is a group action, there isn't an isosceles triangle in any boolean geometry.The geometric results of commutative l-groups \cite{Blu1}, boolean l-groups\cite{Blu2}, and blumenthal\cite{Blu1},\cite{Blu2}, DRl-semigroup\cite{S2} and Ellis \cite{SW A}, \cite{E2} for boolean algebra are generalized by the geometric results found in this study. Blumenthal\cite{Blu1} and penning \cite{pen} studied the metric betweeness in Boolean metric space and swamy\cite{S2} studied the metric betweeness in l-groups. Subbarao\cite{BV} studied the geometry of representable autometrized algebra. Thus, Subbarao generalized the geometric results of swamy \cite{S2} for commutative l-groups and those Blumenthatal and Ellis\cite{E2} for Boolean algebras. In this paper, we introduce the concept of metric betweeness in AL-monoid and discussed the geometric results. And, we obtain several results regarding the properties of metric betweeness, B-linearity and D-linearity in AL-monoid. Also, we obtain some interesting geometric results of DRl-semigroup\cite{S2}. We also introduce the concept of triangle in AL-monoids and obtained the interesting result there do not exist equilateral triangles in AL-monoids. Although a DRl-semigroup is thought of as a vast generalization of commutative l-groups and Brouwerian algebras, we do not obtain many of their geometries for a general DRl-semigroup. In spite of that Swamy studied the geometry of direct product of Boolean algebra and commutative l-groups with some conditions called Boolean l-group.  In this paper we study the Geometry of AL-monoids, by discussing the notion of Metric betweeness, lattice betweeness, B-linearity and D-linearity, segments and equilateral triangles. 
In section \ref{1}, We racall some preliminaries used for our result, and in section \ref{2}, We introduce generalization of DRl-semigroup called Autometrized algebra and  study the Geometry of Autometrized lattice ordered monoids, which subsumes all the geometric properties of DRL-semigroups \cite{S2}. Notation, terminology and results of the paper \cite{T} are employed in this paper. 
\section{preliminaries}	\label{1}
 Here after,throughout this section , $A=(A,+,\leq,\ast,0) $ stands for AL-monoids.In order to commence discussing the geometry of A, we will first introduce the terminology provided in \cite{N} following 
 \begin{definition}\cite{SKLN}
	An Autometrized algebra A is a system $(A,+,\leq,\ast )$ where 
	\begin{enumerate} \item $(A,+)$ is a binary commutative algebra with element 0,
		\item $\leq$ is antisymmetric, reflexive ordering on A,
		\item $\ast : A\times A\rightarrow A$ is a mapping satisfying the formal properties of distance, namely,
		\begin{enumerate}
			\item $a\ast b\geq 0$ for all $a,b$ in A,equality,if and only if $a=b$,
			\item$a\ast b=b\ast a$ for all $a,b$ in A, and
			\item $a\ast b\leq a\ast c +c\ast b$ for all $a,b,c$ in A.
		\end{enumerate}
	\end{enumerate}
\end{definition}
\begin{definition}
	\cite{S2}| A system $A=(A,+,\leq,\ast )$ of arity $(2,2,2)$ is called a Lattice ordered autometrized algebra, if and only if, A satisfies the following conditions.
	\begin{enumerate}
		\item $(A,+,\leq)$ is a commutative lattice ordered semi-group with $'0'$, and
		\item $\ast$ is a metric operation on A.i.e, $\ast$ is a mapping from $A\times A$ into A satisfying the formal properties of distance, namely,
		\begin{enumerate}
			\item $a\ast b\geq 0$ for all $a,b$ in A,equality,if and only if $a=b$,
			\item$a\ast b=b\ast a$ for all $a,b$ in A, and
			\item $a\ast b\leq a\ast c +c\ast b$ for all $a,b,c$ in A.
		\end{enumerate}
	\end{enumerate}
\end{definition}
\begin{definition}
	\cite{BV} A lattice ordered autometrized algebra $A=(A,+,\leq,\ast )$of arity $(2,2,2)$ is called representable autometrized algebra, if and only if, A satisfies the following conditions:
	\begin{enumerate}
		\item $A=(A,+,\leq,\ast)$ is semiregular autometrized algebra. Which means $a\in A$ and $a\geq 0$ implies $a\ast 0=a,$ and 
		\item  for every a in A, all the mappings $x\mapsto a+x,x\mapsto a\vee x,x\mapsto a\wedge x$ and $x\mapsto a\ast x$ are contractions (i.e.,if $\theta$ denotes any one of the operations $+,\wedge,\vee$and $\ast$,then,for each a in A,$(a\theta x)\ast (a\theta y)\leq x\ast y$ for all $x,y$ in A)
	\end{enumerate}
\end{definition}
\begin{definition}\cite{T}\label{2.1} An Autometrized lattice ordered monoid (AL-monoids, for short) is an algebra $(A,+,\vee,\wedge,\ast,0)$ of arity $(2,2,2,2,0)$ where 
	\begin{enumerate}
		\item $(A,+,\vee,\wedge,0)$ is a commutative lattice ordered monoid.
		\item  $a\ast(a\wedge b)+ b  = a\vee b.$
		\item The mappings $x\mapsto a+x,a\vee x, a\wedge x, a\ast x $ are contractions with respect to $\ast$ (A mapping $f:A \longrightarrow A $ is called  a contraction with respect $\ast \Leftrightarrow f(x)\ast f(y) \leq x\ast y $ where $\leq $ is ordering in $A$ induced by $(A,\vee,\wedge)).$
		\item $[a\ast (a\vee b)]\wedge [b\ast (a\vee b)] = 0.$
	\end{enumerate}
 \end{definition}
 \begin{remark}
	Rest of the paper,we simply write $A$ for an AL-monoid  $A=(A,+,\leq,\ast,0)$ and $a,b,c,x,y,z$ stand for elements of $A$.
\end{remark}
\begin{example}
	Any DRl-semigroup is an AL-monoid if $a\ast b=(a-b)\vee (b-a)$.
\end{example}
However, the converse fails. For instance, consider the following
\begin{example}\label{Ex}
	Let $A=\mathbb{Z}\cup \lbrace u\rbrace $, where $\mathbb{Z}$ is the set of all integers and $u$ is an element which is not in $\mathbb{Z}$. For all $a,b\in \mathbb{Z}$ define $+,\ast$ in $A$ as follows: \begin{eqnarray*} a+b = \textrm{the usual sum},\\ a+u=u=u+a,u+u=u,\\a\ast b=|a-b|,a\ast u=u=u\ast a.\end{eqnarray*} Define $\leq $ in $A$ for all $a \in \mathbb{Z}$ as the usual ordering.Then it can be verified that $A=(A,\leq, +,\ast)$ is an AL-monoid which is not a  DRl-semigroup, since there is no least element in $A$ such that $x+u\geq a$ for any $a$ in $\mathbb{Z}.$
\end{example}
The following examples shows that the condition (2) and condition (4) of Definition (\ref{2.1}) are independent in any AL-monoid.
\begin{example}\label{ex}
	Let A be the lattice of all closed subset of the space of real numbers with the usual topology. Then, it can be verified that A satisfies $a\ast (a\wedge b )+(a\wedge b)=a$. but does not satisfy (4).Consider  $a\ast (a\vee b)\wedge (b\ast (a\vee b))=0.$ For, let $x$ and $y$ denote the closed intervals $[0,2]$ and $[2,3]$ respectively. Then, $(x\ast (x\vee y))\wedge (y\ast (x\vee y))=\lbrace 2 \rbrace \neq 0.$ 
\end{example}
\begin{remark}
    From Example \ref{ex} we observe that Representable autometrized algebra is a wider class of AL-monoid and AL-monoid is a sub-class of Representable Autometrized algebra.
\end{remark}
\begin{example}
	Let $\mathbb{Z}\vee \lbrace u,v\rbrace$, where $\mathbb{Z}$ is the set of all intigers and $u$ and $v$ are elements which are not in $\mathbb{Z}:$ Define $+$ in A as follows: 
 \begin{eqnarray*} a+b \textrm{ as usual sum of $a$ and $b$, for all $a,b\in \mathbb{Z}$},\\ a+u=u=u+a,a+v=v=v+a \textrm{for all $a\in \mathbb{Z}$}, u+v=u=v+u, u+u=u,v+v=v \end{eqnarray*}. Define $\leq $ in $A$ as follows: for the elements in $\mathbb{Z}$, let $\leq $ be the usual ordering, define $u<a<v, \forall a \in \mathbb{Z}.$ Define $\ast $ in $A$ as follows:\\ $a\ast b=|a-b|,\forall a,b\in Z, a+u=v=v+a \forall a\in \mathbb{Z}, a+v=v=v+a$ for all $a$ in $\mathbb{Z}, u\ast v=v=v\ast u$ and $u\ast y=0=v\ast v.$ Then, it can be verified that $(A,+,\leq, \ast )$ satisfies axioms $(1),(3),(4),$ of Definition \ref{2.1} and but does not satisfy axiom $(2)$, because $v\ast (v\wedge u)+v\wedge u=v\ast u+u=v+u=u\neq v.$
\end{example}
The following theorem shows that any AL-monoid can be equationally definable  so that the class of AL-monoids are closed under the formation of subalgebras, direct unions, and homomorphic images and hence form a variety.
\begin{lemma}\label{l1}
$a\leq b$ implies $a\ast c\leq b\ast c,$ for any $ c\in A.$
	\begin{proof}
From $b\ast (b\wedge c)+c= b\vee c$ we get  $b\ast c+c=b,$ and $a\ast c+c=a.$ Hence, $a\leq b \Rightarrow a\ast c +c \leq b\ast c +c$ implies $a\ast c\leq b\ast c$(since $(A,+,\leq )$ is commutative l-group.  
	\end{proof}
\end{lemma}
\section{Result}\label{2}
\begin{definition}
	A subset $S$ of an AL-monoid $A =(A, + , \leq,\ast,0 )$ is said to be a $~'$subgeometry' if and only if $a$ and $b$ belong to $S$ implies that $a\ast b$ belongs to S.
\end{definition}
We shall treat any distinct elements $a,b,c$ of an AL-monoid $(A,+,\leq,\ast,0)$ as the vertices of triangle with sides $a\ast b, b\ast c$ and $c\ast a.$\\
\begin{theorem}
	For $\triangle(a,b,c)$ in $A$, the following are equivalent:
	\begin{enumerate}
		\item[i] $\triangle(a,b,c)$ has fixty.
		\item[ii] $a\vee b=b\vee c=c\vee a$ and $a\wedge b\wedge c=0.$
		\item[iii] $\triangle(a\ast b,b\ast c,c\ast a)$ has fixty and the set $\lbrace 0,a,b,c\rbrace$ is a subgeometry where (where $x\ast y=|x-y|).$
	\end{enumerate}
	\begin{proof}
$(i)\implies (ii):$\\
 Assume that $\triangle (a,b,c)$ has fixty. So, we have $a\ast b=c,b\ast c=a,$ and $c\ast a=b.$ Clearly, $a>0,b>0$ and $c>0$. Now, $a\vee b=a\ast b+a\wedge b=c+a\wedge b\geq c,$ since $a\wedge b\geq 0.$\\
 Similarly ,$b\vee c\geq a$ and $c\vee a \geq b$. Hence, $a\vee b=b\vee c=c\vee a.$ \\ Now, $$a\vee b=a\ast b+a\wedge b=c+a\wedge b\geq c,$$ since ($a\wedge b\geq 0)$. \begin{eqnarray*}
     b\vee c&=& b\ast c +b\wedge c \\&=& a+ b\wedge c \\&\geq a&
 \end{eqnarray*} and \begin{eqnarray*}
     c\vee a&=&c\ast a +c\wedge a \\&=& b+c\wedge a \\&\geq& b. 
 \end{eqnarray*} Hence, $$a\vee b=b\vee c= c\vee a.$$
 Now, $$a=b\ast c\leq b\ast (a\vee b)+(a\vee b)\ast c =b\ast (a\vee b)+(a\vee c)\ast c$$,and $$c=a\ast b\leq a\ast (a\vee c)\ast a.$$ Hence, $$
 a\wedge c \leq (b\ast (a\vee b))+(a\vee c)\ast c) \wedge (b\ast (a\vee b))+(a\vee c)\ast a$$$$= b\ast (a\vee b)+((a\vee c)\ast c)\wedge ((a\wedge c)\ast a)= b\ast (a\vee b).$$\\
Further, $$b=a\ast c\leq a\ast (a\vee b)+(a\vee b)\ast c$$$$=a\ast (a\vee b)+(b\vee c)\ast c$$, and $$c=a\ast b\leq (a\ast (a\vee b)+(a\vee b)\ast b)$$$$=(a\ast (a\vee b)+(a\vee b))+((b\vee c)\ast b).$$ \\ Hence, $$b\wedge c\leq (a\ast (a\vee b)+(b\vee c)\ast b)$$$$=(a\ast (a\vee b)+(b\vee c)\ast c)\wedge (b\vee c)\ast b=a\ast (a\vee b).$$  Hence, $a\wedge b\wedge c \geq 0,$ it follows that  $a\wedge b\wedge c=0.$ Thus, $(i)\implies (ii).$\\
\newline	$(ii)\implies (i).$ Assume that $a\vee b=b\vee c=c\vee a$ and $a\wedge b\wedge c=0.$ clearly, $a\geq 0,b\geq 0$ and $c\geq 0.$ So, we have \begin{eqnarray*} c&=&c\ast 0\\&=&c\ast (c\wedge a\wedge b)\\&=&(c\wedge (a\vee c))\ast (c\wedge a\wedge b)\\ &\leq& (a\vee c)\ast (a\wedge b) \\&=&(a\vee b)\ast (a\wedge b)\\ &=&a\ast b\end{eqnarray*}. Thus, $c\leq (a\ast b)$. Now,\begin{eqnarray*} a\ast b&\leq& a\ast (a\vee c)+(a\vee c)\ast b\\&=&a\ast (a\vee b)+(b\vee c)\ast (b\vee 0)\\&\leq& a\ast (a\vee b)+(c\ast 0)\\&=& a\ast (a\vee b)+c.\end{eqnarray*} Hence, \begin{eqnarray*} a\ast b&\leq& (c+a\vee b)\ast b)\wedge (c+a\ast (a\vee b))\\ &=& c+((a\vee b)\ast b)\wedge ((a\vee b)\ast a)\\&=&0 \end{eqnarray*}.Thus,$a\ast b\leq c$. Hence, $\triangle(a,b,c)$ has fixty. Thus, $(2)\implies (1).$Hence, $(1)$ and $(2)$ are equivalent.\\
$(i)\implies (iii)$. $\triangle(a,b,c)$ has fixty implies $a\ast b=c,b\ast c=a, a\ast c=b.$ This implies $$(a\ast b)\ast (b\ast c)=(a\ast (b\ast b)=a\ast 0\ast c=a\ast c,$$$$(a\ast b)\ast (c\ast a)=(a\ast a)\ast (b\ast c)=0\ast (b\ast c)=b\ast c,$$ $$(b\ast c)\ast (c\ast a)=b\ast (c\ast c)\ast a=b\ast 0\ast a=b\ast a.$$
To show $S=\lbrace 0,a,b,c\rbrace$ is a subgeometry of $A$,since from (i) $\triangle(a,b,c)$ has fixty $a\ast b=c\in S, a\ast c=b\in S,b\ast c=a, 0\ast a=a\in S,0\ast b=c\in S,0\ast c=c\in S,$ shows that S is a subgeometry.
	\newline $(iii)\implies (i)$: assume that $\triangle (a\ast b, b\ast c, c\ast a)$ has fixty and $\lbrace 0,a,b,c\rbrace$ is a subgeometry. Now, $a\ast b$ should by any of $a,b,c$(since $a\neq b\neq c\neq a$).\\ Let $a\ast b=a$. If $b\ast c=a,$ then, $a\ast c=(a\ast b)\ast (b\ast c)=a\ast a=0,$ so that $a=c$.\\ So, $b\ast c\neq a.$ Hence $b\ast c=b$ or $c$. \\ If $b\ast c=b,$ then, $a\ast c=(a\ast b)\ast (b\ast c)=a\ast b=a,$ so that $b\ast c=(a\ast b)\ast (a\ast c)=a\ast a=0,$ which implies $b=c.$ Hence, $b\ast c\neq b.$ So, $b\ast c=c.$ Then, $a\ast c$ should be any of $a,b,c.$ If $a\ast c=a,$ then, $c=b\ast c$, then, $c=b\ast c=(a\ast b)\ast (a\ast c)=a\ast a=0,$ which implies $b=c=0.$ So, $a\ast c\neq a.$ \\If $a\ast c=b,$ then $a=a\ast b=(a\ast c)\ast (c\ast b)=b\ast c=c.$ So,$a\ast c\neq b.$ If $a\ast c= c,$ then, $a=a\ast b=(a\ast c)\ast (c\ast b)=c\ast c=0,$ which implies $a=b=0.$ Thus, $b\ast c\neq c.$ Hence, $a\ast b\neq a.$ \\Similarly, we can show that $a\ast b\neq b.$ Hence, $a\ast b=c.$ Similarly, $b\ast c=a,$ and $c\ast a=b$. Hence, $\triangle(a,b,c)$ has fixty. Thus, $(iii)\implies (i).$ Hence, $(i)$ and $(iii)$ are equivalent.\\
		 These complete the proof.
	\end{proof}
\end{theorem}
In any Boolean geometry, since the metric operation is group operation , there do not exist any isosceles triangles. In commutative l-groups also, Swamy \cite{S2} has shown that equilateral triangle do not exist.
\begin{theorem}
	There do not exist equilateral triangle in AL-monoid $A.$
	\begin{proof}
		Assume that $a\ast b=b\ast c=c\ast a.$ Now, $c=c\ast (c\wedge a)+c\wedge a=a\ast b+a,$ and $c=c\ast (c\wedge b)+c\wedge b\leq c\ast b+b=a\ast b+b.$ 
		Hence, $c\leq ((a\ast b)+a)\wedge (a\ast b+b)=a\ast b+a\ast b= a\vee b.$ Thus, $c\leq a\vee b$. Similarly, $b\leq a\vee c$ and $a\leq b\vee c$.
		Hence, $a\vee b=b\vee c=c\vee a.\\
		$Now, \begin{eqnarray*} a\ast b & \leq & a\ast (a\vee c)+(a\vee c)\ast b \\ &=& a\ast (a\vee c)+(a\vee b)\ast b\\ &=&a\ast (a\vee c)+(a\vee b)\ast b,\end{eqnarray*} and \begin{eqnarray*} a\ast b&=&b\ast c\\ &\leq & (b\ast (a\vee b))+(a\vee b)\ast c\\ &=& b\ast (a\vee b)+(a\vee c)\ast c.\end{eqnarray*} Hence, $ a\ast b
  \leq(a\ast (a\vee c)+(a\vee b)\ast b) \wedge (b\ast (a\vee b))+(a\vee c)\ast c= (a\ast(a\vee c))\wedge ((a\vee c)\ast c)+(a\vee b)\ast b=(a\vee b)\ast b.$ \\
  Similarly, we can prove that $a\ast b\leq (a\vee b)\ast a.$ Hence, $a\ast b\leq ((a\vee b)\ast b)\wedge ((a\vee b)\ast a)=0.$ Hence, $a\ast b=b\ast c=c\ast a=0,$ which implies $a=b=c.$ Hence, the $\triangle(a,b,c)$ degenerates (i.e., $a,b,c$ are not distinct). This completes the proof.
	\end{proof} 
\end{theorem}
\begin{theorem}
	$A$ is a chain if and only if it is free of triangles with fixty.
	\begin{proof}
		Assume that $A$ is not a chain. So, there exist $a,b \in A$ such that $a\neq b, a\vee b\neq a,$ and $a\vee b\neq b.$ Then, $(a\vee b)\ast a, (a\vee b)\ast b,$ and $a\ast b$ are three distinct elements in $A$.\\ For, if $(a\vee b)\ast a=(a\vee b)\ast b,$ then, $a\ast b\leq a\ast (a\vee b)+(a\vee b)\ast b=2(a\ast (a\vee b))=2(b\ast (a\vee b)),$ and hence, $a\ast b\leq (2(a\ast(a\vee b)))\wedge (2(b\ast (a\vee b)))=0.$ (since in commutative lattice oredered semigroup with identity 0, $x\wedge y=0$ implies $mx\wedge ny=0$ for all positive intigers $m$ and $n$), which implies $a=b$, a contradiction.\\
		\newline
		 If $a\ast b=(a\vee b)\ast a,$ then, $(a\vee b)\ast b\leq a\vee b=(a\vee b)\ast a,$ which implies $(a\vee b)\ast b=((a\vee b)\ast b)\wedge ((a\vee b)\ast a)=0$, so that $a\vee b=b,$ a contradiction. \\ Similarly, if $a\ast b=(a\vee b)\ast b,$ we leads to a contradiction that $a\vee b=a$. Hence, the three elements $a\ast (a\vee b),b\ast (a\vee b),$and hence $a\ast b$ are distinct and so, they form  triangle $\triangle_{1}=\triangle(a\ast (a\vee b)),b\ast (a\vee b),(a\ast b).$\\
		 \newline
		  Now,we show that $\triangle_{1}$ has fixty. For, $(a\ast (a\vee b))\ast (b\ast (a\vee b))\leq a\ast b,$ and $a\ast b=(a\ast (a\wedge b)+a\wedge b)\ast (b\ast (a\wedge b)+a\wedge b)\leq (a\ast (a\wedge b))\ast (b\ast (a\ast (a\vee b))),$ which implies  $(a\ast (a\vee b))\ast (b\ast (a\vee b))=a\ast b.$ Also, $(a\ast (a\vee b))\ast (a\ast b)\leq (a\vee b)\ast b,$ and \begin{eqnarray*}(a\vee b)\ast b&=&(a\ast b+a\wedge b)\ast (a\wedge b)+a\wedge b)\\ &\leq & (a\ast b)\ast (b\ast (a\wedge b))\\ &=& (a\ast b)\ast (a\ast (a\vee b))\\ &=& (a\ast (a\vee b))\ast (a\ast b),\end{eqnarray*} which imply $(a\ast (a\vee b))\ast (a\ast b)=b\ast (a\vee b).$ Similarly, we can prove that $(b\ast (a\vee b))\ast (a\ast b)=a\ast (a\vee b).$ Thus, $\triangle_{1}$ has fixty. So, $A$ is free of triagles with fixty. That is, if $A$ is free of traigles with fixty, then, $A$ is chain. Conversely, if $A$ is chain, there is no traigle in $A$ have fixty.
	\end{proof}
\end{theorem}

\begin{theorem}
	If all traigles in $A$ are isosceles, then, $A$ is chain. 
	\begin{proof}
		For any $a,b$ in $A$, consider the triangle $\triangle(a\vee b,a\wedge b,a)$, now,$(a\vee b)\ast (a\wedge b)=(a\vee b)\ast a$ implies \begin{eqnarray*}a\vee b&=& a\ast b+a\wedge b\\ &=& (a\vee b)\ast (a\wedge b)+(a\wedge b)\\ &=& (a\vee b)\ast a+a\wedge b\\ &=& b\ast (a\wedge b)+a\wedge b\\ &=&b.\end{eqnarray*} Also, $(a\vee b)\ast (a\wedge b)=(a\wedge b)\ast a$ implies $a\vee b=a\ast b+a\wedge b=a\ast (a\wedge b)+(a\wedge b)=a.$ And, $(a\vee b)\ast a=(a\wedge b)\ast a$ implies \begin{eqnarray*} b&=& b\ast (b\wedge a)+b\wedge a\\ &=& a\ast (a\wedge b)+b\wedge a\\ &=& a\ast (a\wedge b)+(a\wedge b)\\ &=& a\end{eqnarray*}. Hence, the triangle $\triangle (a\vee b,a\wedge b,a)$ degenerates(i.e., $a\vee b$,$a\wedge b$, a are not distincts). Hence, $A$ is a chain.
	\end{proof}
\end{theorem}
\begin{definition}
	An element $x\in A$ said to $''$ lie metrically between $'$a' and $'$ b' (in symbols $(a,x,b)M)$ if and only if $a\ast x+x\ast b=a\ast b.$
\end{definition}
\begin{definition}
    Any betweeness relation $B$ is said to satisfy the special inner property $\beta: B(a,b,c)$ and $B(a,c,b),$ if and only if, $b=c.$
\end{definition}
\begin{lemma}
	In $A$, three points $a,b,c$ fail to have special inner property, if $a,b,c$ are the vertices of an isosceles triangle in which the sum of any two sides is equal to their union(join).
	\begin{proof}
		Assume that $a,b,c$ are the vertices of an isosceles traingle, for which $(a\ast b)\vee (a\ast c)=a\ast b+a\ast c=a\ast b+a\ast c, (a\ast b)\vee (b\ast c)=(a\ast b)+(b\ast c), (a\ast c)\vee (b\ast c)=a\ast c+b\ast c,$ and $a\ast b=b\ast c.$ Thus, $a\ast b+b\ast c=(a\ast b)\vee (b\ast c)\leq (a\ast b)\vee (b\ast a+a\ast c)$(by $(M_{3}))$. SO, $a\ast c+c\ast b=a\ast b+c\ast b=a\ast c=a\ast b.$ Thus, we have $(a,b,c)M$ and $(a,c,b)M.$ But, $b\neq c$,since $a,b,c$ are the vertices of triangle, the special inner property fails in this case. This completes the proof.
	\end{proof}
\end{lemma}
\begin{definition}
    In AL-monoid $A,$ any betweeness relation $M$ is said to have 
    (i). transitive $t_{1},$ if and only if, $B(a,b,c), B(a,d,b)\implies B(d,b,c).$, (ii). the transitive $t_{2},$ if and only if, $B(a,b,c), B(a,d,b)\implies B(a,d,c).$
\end{definition}
\begin{theorem}
	If the metric betweeness has transitivity $t_{1}$, then, $A$ is free of isosceles triangles for which the sum of any two sides is equal to their union.
	\begin{proof}
		$(a,b,c)M$ and $(a,d,b)M$ imply $a\ast b+b\ast c=a\ast c$ and $a\ast d+d\ast b=a\ast b,$ which imply $a\ast d+d\ast b+ b\ast c=a\ast b+b\ast c=a\ast c,$ so that $a\ast c\geq a\ast d+d\ast c.$ But, $a\ast c\leq a\ast d+d\ast c$.
	\end{proof}
\end{theorem}
\begin{lemma}
	In $A$, the metric betweeness has transitivity $t_{2}$.
	\begin{proof}
		$(a,b,c)M$ and $(a,d,b)M$ imply $a\ast b+b\ast c=a\ast c$ and $a\ast d+d\ast b=a\ast b,$ which imply $a\ast d+d\ast b+b\ast c=a\ast c,$ so that $a\ast c\geq a\ast d+d\ast c.$ But, $a\ast c\leq a\ast d+d\ast c).$ Hence, $(a,d,c)M$
		This completes the proof.
	\end{proof}
\end{lemma}
\begin{theorem}\label{T12}
	In $A$, lattice betweeness implies metric betweeness.
\end{theorem}
\begin{lemma}
	$a\ast b+b\ast c=(a\wedge c)\ast b+b\ast (a\vee c)$ for all $a,b,c$ in $A$.
	\begin{proof}
		$a\ast b+b\ast c\leq a\ast (a\vee c)+(a\vee c)\ast b+b\ast (a\wedge c)+(a\wedge c)\ast c=((a\vee c)\ast b+b\ast (a\wedge c))\ast (a\ast (a\vee c)+a\ast (a\vee c))=((a\vee c)\ast b\ast b\ast (a\wedge c))+2(a\ast(a\vee c))$ and similarly, $a\ast b+b\ast c \leq ((a\vee c)\ast b+b\ast (a\wedge c))+2(c\ast (a\vee c))$ and similarly, $a\ast b+b\ast c\leq ((a\vee c))\leq ((a\vee c)\ast b+b\ast (a\wedge c))+2(c\ast (a\vee c)).$\\ 
		\newline Hence, $a\ast b+b\ast c\leq ((a\vee c)\ast b+b\ast (a\wedge c))+2(a\ast (a\vee c))\wedge ((a\vee c)\ast b+b\ast (a\wedge c))+2(c\ast (a\vee c))=((a\vee c)\ast b+b\ast (a\vee c))$. Since in any commutative lattice ordered semigroup with identity $a\wedge x=0$ and $a\wedge y=0$ imply $a\wedge (x+y)=0$. Further, $(a\vee c)\ast b+b\ast (a\wedge c)\leq (a\vee c)\ast a +a\ast b+b\ast c+c\ast (a\wedge c)=((a\ast b)+b\ast c)+2(a\ast (a\vee c))$ and similarly $(a\vee c)\ast b +b\ast (a\wedge c)\leq (a\ast b+b\ast c)=2((c\ast (a\vee c))).$ \\ Hence, we have $a\ast b+b\ast c=(a\vee c)\ast b+b\ast (a\wedge c)=(a\wedge c)\ast b+b\ast (a\vee c)$ for all $a,b,c$ in $A.$ This completes the proof.
	\end{proof}
\end{lemma}
\begin{corollary}
	In $A$, $(a,b,c)M$ iff $(a\wedge c, b, a\vee c)M.$
\end{corollary}
\begin{lemma}
	In $A$,$a\leq b\leq c$ implies $(a,b,c)M.$
	\begin{proof}
		Assume that $a\leq b\leq c.$ Then, $c\ast b=(c\ast (c\wedge a))+(c\wedge a)\ast (b\ast (b\wedge a)+b\wedge a)=(c\ast a+a)\ast (b\ast a+a)\leq (c\ast a)\ast (b\ast a)\leq c\ast b.$ Hence, $c\ast b=(c\ast a)\ast (b\ast a).$ And,$b\ast a\leq c\ast a,$ since, $b\ast a=(b\wedge c)\ast (a\wedge b)\leq c\ast a.$ So, $c\ast b+b\ast a=(c\ast a)\ast (b\ast a)+b\ast a=(c\ast a)\ast ((c\ast a)\wedge (b\ast a))(c\ast a)\wedge (b\ast a)=c\ast a.$ This implies $(a,b,c)M.$ This complete the proof.
	\end{proof}
\end{lemma}
\textbf{proof of the Theorem \ref{T12}}\\
Assume that $(a,b,c)L$. So, $a\wedge c\leq b\leq b\vee c.$ Hence, $(a\wedge c)\ast b+b\ast (a\vee c)=a\ast c.$ Hence $(a,b,c)M$. This complete the proof.\\
lattice betweeness and metric betweeness are not equivalent in $A$. However, we have the following
\begin{theorem}
	In $A$, lattice betweeness and metric betweeness are equivalent iff metric betweeness has transitivity $t_{1}$.
	\begin{proof}
		Assume that metric betweeness has transitivity $t_{1}$, it is enough if we verify $(1)$ through $v$ of metric betweeness. Clearly, metric betweeness satsisfies (i) and $(2)$. Since the metric betweeness has transitivity $t_{1}$, $(ii)$ holds, because, $(a,b,c)M$ and $(a,c,b)$M, i.e., $c=b$. by above lemma $(iii)$ is also satisfied. Since we have $(a\wedge c, a\vee c, a\vee c)M$ and $(a\wedge c, a\wedge c, a\vee c)$, by corollary, (4) is satisfied. by lemma() $(5)$ is satisfied. Hnece the metric betweeness and lattice betweeness are equivalent.\\
		Conversely, assume that lattice betweeness and metric betweeness are equivalent. Since trhe lattice betweeness has transitivity $t_{1}$ in any lattice. Follows that metric betweeness has transitivity $t_{1}$. This completes the proof.
	\end{proof}
\end{theorem}

\begin{definition}
	AL-monoid $A=(A,+,\leq, \ast )$is said to be "metrically convex" if and only if given $a\neq b$ in $A$ such that $a\neq x\neq b$ and $a\ast x+x\ast b=a\ast b.$ And, we say that an element $c$ in $A(c\neq 0)$ is an atom if and only f $0<x<c$ is imposible for any $x$ in $A$
\end{definition}
\begin{definition}
    An AL-monoid $A$ is ptolemaic, if and only if, for every four points $a,b,c,d$ in $A,$ it is true that 
    \begin{eqnarray*}
        (a\ast b)\wedge (c\ast d)\leq a\ast c \wedge b\ast d +(a\ast d)\wedge (b\ast c\\ (a\ast c)\wedge (b\ast d)\leq (a\ast b)\wedge (c\ast d)+(a\ast d)\wedge (c\ast b)\\ (a\ast d)\wedge (b\ast c)\leq (a\ast c)\wedge c \wedge d\ast b+a\ast b\wedge c\ast d.
    \end{eqnarray*}
\end{definition}
\begin{theorem}
	Al-monoids are ptolemaic with $\wedge$ as multiplication.
\begin{proof}
Let $a,b,c,d \in A.$ Let $a=a_{1}\ast a_{2}, b=a_{2}\ast a_{3}, c=a_{3}\ast a_{1}, x=a_{3}\ast a_{4},y=a_{4}\ast a_{1},$ and $z=a_{2}\ast a_{4}.$ We have , $a\leq b\ast c$ and $x\leq b+z$ and hence $a\wedge x \leq (b+c)\wedge (b+z)=b+(c\wedge z)$. Also, $a\leq a+z$ and $x\leq c+y\implies a\wedge x\leq (y+c)\wedge (y+z)=y+(c\wedge z)$.$ a\wedge x\leq b+(c\wedge z)$ and $a\wedge x\leq y+c\wedge z\implies a\wedge x\leq (b+c\wedge z)\wedge (y+c\wedge z)=b\wedge y +(c\wedge z).$Also, $(a\ast c)\wedge (b\ast d) +a\ast d \wedge b\ast c=(a\ast c)+((a\ast d)\wedge (b\ast c))\wedge (b\ast d +(a\ast d\wedge (b\ast c)))=((a\ast c)+(a\ast d))\wedge ((a\ast c +b\ast c))\wedge ((b\ast d)+(a\ast d))\wedge (b\ast d+b\ast c))\geq c\ast d \wedge (a\ast d)\wedge (b\ast a)\wedge (d\ast c)\geq c\ast d \wedge a\ast b=(a\ast b)\wedge (c\ast d)$. We can proof the other in similar ways. Since, in $A,$ we can not find pair wise disjoint elements $a,b,c$ in $A$ such that $a\ast b=b\ast c=c\ast a,$ there do not exist equilateral triangle in $A.$
\end{proof}
\end{theorem}
\begin{definition}
An n-tuple of distinct points in AL-monoids $A$ for $(n\geq 3)$ is said to be 
    \begin{enumerate}
        \item[i.] B-linear, if and only if, there exists a labelling $(P_{1},p_{2},...,P_{n})$ of its elements such that $B(P_{i},P_{j},P_{k})$ holds wherever $1\leq i<j<j<k\leq n).$
        \item[ii.] D-linear, if and only if, there exists a labeling $(P_{1},p_{2},...,P_{n})$  of its elements such that $P_{1}\ast P_{n}=P_{1}\ast P_{2} +P_{2}\ast P_{2}\ast P_{3}+...+p_{n-1}\ast P_{n}.$
    \end{enumerate}
\end{definition}
\begin{theorem}
		In $A,$ B-linearity implies D-linearity.
	\begin{proof}
		Assumet that an n-tuple $p$ is a Bi-linear, then, there exists a labelling $(p_{1},p_{2},...,p_{n})$ of elements of $p$ such that $B(p_{i},p_{j},p_{k})$ holds whenever $i\leq i<j<k\leq n$, from which  follows, in particular, that $B(p_{i},p_{i+1},p_{n})(i<n-1)$ holds, so that $B(p_{i},p_{j},p_{k})M$ by Theorem 3.14, which implies $P_{i}\ast p_{n}=p_{i}\ast p_{i+1}+p_{i+1}\ast p{n}$ for $i<n-1.$ Hence, $p_{1}\ast p_{n}=p{1}\ast p_{2}+p_{2}\ast p_{n}=p_{1}\ast p_{2}\ast p_{3}\ast p-{n}=\sum_{i=1}^{n-1}(p_{i}\ast p_{i+1}),$ so that $P$ is D-linear.This complete proof.
	\end{proof}
		\end{theorem}
\begin{theorem}
	In $A,$ the following assertions are equivalent:
	\begin{enumerate}
		\item $(a,b,c)L$ if and only if $(a,b,c)M$.
		\item Metric betweeness has transitivity $t_{1}.$
		\item $a\leq b\vee c$ and $a\ast b\geq a\ast c$ imply $b\leq c.$
		\item $(a,b,c)M$ and $(a,c,b)M$ imply $b=c.$
\begin{proof}
$(1)\Longleftrightarrow (2)$, by Theorem (3.18). Now, we prove $(1)\implies (4)\implies (3)\implies (1,)$\\$(1)\implies (4)$ Assume that $(a,b,c)L$ iff $(a,b,c)M$. Suppose $(a,b,c)M$ and $(a,c,b)M.$ Then, by assumption, we have $(a,b,c)L$ and $(a\wedge c\leq c\leq a\vee b).$ So, $a\wedge c\leq a\wedge b\leq a\wedge c,$ and $a\vee c\leq a\vee c.$ Hence, $a\wedge c=a\wedge b,$ and $a\vee c=a\vee b.$ Since $(A,\leq)$ is distributive lattice (by Theorem ), it follows that $b=c$. Thus, $(1)\implies (4)$.\\
	\newline
$(4)\implies (3)$. Assume that $(a,b,c)M$ and $(a,c,b)M$ imply $b=c.$Suppose $a\geq b\vee c$ and $a\ast b\geq a\ast c.$ First, we prove that $a\ast (b\wedge c)=a\ast b$.For, $a\ast (b\wedge c)\leq a\ast c+c\ast (b\wedge c)($by $M_{3})\leq a\ast b +c\ast (b\wedge c)$. So, $a\ast (b\wedge c)\leq (a\ast b)+b\ast (b\wedge c)\wedge(a\ast b+c\ast (b\wedge c)) =a\ast b+(c\ast (b\vee c))\wedge (b\ast (b\vee c))=a\ast b.$ \\And, $a\ast b=(a\vee b)\ast (b\vee (b\wedge c))$(by $L_{2})$. Hence, $a\ast (b\wedge c)=a\ast b.$ Since $b\wedge c\leq b\leq a,$ by lemma 3.17, we have $(b\wedge c,b,a)M$, i.e., $(a,b,b\wedge c)M$. So, $a\ast (b\wedge c)+(b\wedge c)\ast b=a\ast b+b\ast (b\wedge ca\ast (b\wedge c)=a\ast b,$ i.e., $(a.b\wedge c. b)M.$ Thus, we have $(a,b,b\wedge c)M$ and $(a,b\wedge c,b)M.$ Hence, by our assumption, it follows that $b=b\wedge c.$ i.e., $b\leq c$. Thus, $(4)\implies (3)$.\\
	\newline 
$(3)\implies (1)$. Assume that $a\geq (b\vee c)$ and $a\ast b\geq a\ast c$ imply $b\leq c$ (by Theorem (3.14) ). Now, assume $(a,b,c)M$  that is $a\leq c$ ;$b=b\ast (b\wedge a)+b\wedge a\leq b\ast a+a\leq c\ast a.$ (since $(a,b,c)M$ means $a\ast b+b\ast c=a\ast c=c\ast (c\wedge a)+c\wedge a=c.$ Thus, $a\leq c$ and $b\leq c,$ and hence, $c\geq a\vee b.$ But, $c\ast a\geq c\ast b.$ Hence, by our assumption, it follows that $a\leq b.$ Thus, $a\leq b\leq c.$ Hence, $(a,b,c)L.$ Thus, we have proved $(1)\implies (4)\implies (3)\implies (1)\implies (2).$ Hence, $(1),(2),(3)$ and $(4)$ are equivalent.
		\end{proof}
	\end{enumerate}
\end{theorem}
\begin{theorem}
    Representable autometrized algebra is larger class of AL-monoids.
    \begin{proof}
        \begin{example}
     Let $A=[0,2], B=[1,3]$ and $C=[0,1]\vee [2,3]$ then $\triangle(A,B,C)$ has fixty, while $A\wedge B\wedge C\neq \emptyset. $       
        \end{example}
    \end{proof}
\end{theorem}


\end{document}